\RequirePackage{fix-cm}

\documentclass[smallextended]{svjour3}       
\smartqed  

\usepackage{amssymb}
\usepackage{graphicx}
\usepackage{amsmath}
\usepackage{multirow,url}
\usepackage{float}
 \usepackage{rotating}
 \usepackage{mathabx}
\usepackage{mathptmx}      
\usepackage{xcolor}
\usepackage{latexsym}
\usepackage{hyperref}
\hypersetup{hypertex=true,
colorlinks=true,
linkcolor=blue,
anchorcolor=blue,
citecolor=blue}

\usepackage{caption}
\usepackage{subcaption}

\begin{document}
\title{Absolute-value based preconditioner for complex-shifted Laplacian systems}
\author{\small Xuelei Lin\and  Congcong Li \and Sean Hon 
}


\institute{Xuelei Lin \at
              School of Science, Harbin Institute of Technology,	Shenzhen 518055, China\\
       \email{linxuelei@hit.edu.cn}
\and
Corresponding author: Congcong Li \at Department of Mathematics, Hong Kong Baptist University, Kowloon Tong, Hong Kong SAR \\             \email{22482245@life.hkbu.edu.hk}      
\and
          Sean Hon \at
              Department of Mathematics, Hong Kong Baptist University, Kowloon Tong, Hong Kong SAR\\
        \email{seanyshon@hkbu.edu.hk}        
}

\date{Received: date / Accepted: date}

\maketitle

\abstract{The complex-shifted Laplacian systems arising in a wide range of applications.
In this work, we propose an absolute-value based preconditioner for solving the complex-shifted Laplacian system. In our approach, the complex-shifted Laplacian system is equivalently rewritten as a $2\times 2$ block real linear system. With the Toeplitz structure of uniform-grid discretization of the constant-coefficient Laplacian operator, the absolute value of the block real matrix is fast invertible by means of fast sine transforms. For more general coefficient function, we then average the coefficient function and take the absolute value of the averaged matrix as our preconditioner. 
With assumptions on the complex shift, we theoretically prove that the eigenvalues of the preconditioned matrix in absolute value are upper and lower bounded by constants independent of matrix size, indicating a matrix-size independent linear convergence rate of MINRES solver. Interestingly, numerical results show that the proposed preconditioner is still efficient even if the assumptions on the complex shift are not met. The fast invertibility of the proposed preconditioner and the robust convergence rate of the preconditioned MINRES solver lead to a linearithmic (nearly optimal) complexity of the proposed solver. The proposed preconditioner is compared with several state-of-the-art preconditioners via several numerical examples to demonstrate the efficiency of the proposed preconditioner. }




\maketitle

\section{Introduction}\label{sec1}
In this work, we consider iteratively solving the following  complex-shifted Laplacian system:
\begin{equation}\label{eqn:complex_shifted_Lap_system}
  (K + \lambda I)\mathbf{z}=\mathbf{f},
\end{equation}
where $K \in \mathbf{R}^{M \times M}$ is derived from the discretization of the Laplace operator $-\nabla(a(\mathbf{x}) \nabla)$ based on finite difference scheme on rectangular domains and $0 < \check{a}\leq a(x) \leq \hat{a}$, $\lambda=\alpha+\beta {\bf i}$ is a complex number with ${\bf i}$ denoting the imaginary unit and $\alpha,\beta\in\mathbb{R}$, $I$ is identity matrix of proper size, $\mathbf{z}$ is the unknown vector to be solved, and $\mathbf{f}$ is the given right-hand side vector. 

The complex-shifted Laplacian matrix was initially introduced as a ``shifted Laplace" preconditioner for handling the Helmholtz equations as documented in \cite{erlangga2004class}. As a preconditioner, its application follows with the computation of its inverse times some given vectors, from which complex shifted Laplacian system \eqref{sec1} arises. Besides, the complex-shifted Laplacian system \eqref{eqn:complex_shifted_Lap_system} also arises other scenarios, such as preconditioned iterative solutions of local and non-local evolutionary equations with Laplacian operator as spatial term; see, e.g., \cite{he2022vanka,mcdonald2018preconditioning,lin2021all,gander2020paradiag,liu2020fast,lin2018separable,lin2016fast,lu2015fast,liu2022well}. In these scenarios, the triangular matrices arising from discretization of temporal derivatives are further approximated by $\alpha$-circulant matrices which have complex eigenvalues. As a result, the approximated time-space matrices are block diagonalizable with each eigen-block being a complex-shifted Laplacian matrix. Solving linear systems with the approximated time-space coefficient matrices gives rise to complex-shifted Laplacian systems of form \eqref{eqn:complex_shifted_Lap_system}. These applications motivate us to study fast solvers for \eqref{eqn:complex_shifted_Lap_system}.

Although direct methods have the potential to offer benefits, solving complex-shifted linear systems \eqref{eqn:complex_shifted_Lap_system} using these methods can be computationally expensive. This serves as strong motivation for us to explore and develop more efficient solvers for such systems. The efficient computation of Helmholtz-like complex symmetric systems has been extensively explored in the literature. For the multigrid method, notable references include \cite{cocquet2017large,maclachlan2008algebraic,notay2010aggregation,wu2018solving,hocking2021optimal}.  In addition, for the preconditioned Krylov subspace
 methods like GMRES or BiCGStab methods, relevant contributions can be found in \cite{gander2015applying,graham2017domain,erlangga2004robust}. {In \cite{bai2007successive,bai2011preconditioned,axelsson2014comparison}, the Preconditioned Modified Hermitian and Skew-Hermitian Splitting method (PMHSS) has been proposed for the complex linear system \eqref{eqn:complex_shifted_Lap_system}, which is based on Hermitian and Skew-Hermitian splitting of the coefficient matrix. However, the corresponding theoretical analysis relies on assumptions that $K+\alpha I$ is symmetric positive definite (SPD) and that $\beta I$ is symmetric positive semi-definite (SPSD).} Recently, the authors generalized the additive element-wise Vanka smoother \cite{hocking2021optimal} for solving complex-shifted Laplacian systems in ParaDIAG algorithms \cite{he2022vanka}. Nevertheless, the convergence theory for these iterative solvers has considerable room to be developed.

Instead of directly working on \eqref{eqn:complex_shifted_Lap_system}, we study a preconditioner for the following equivalent $2\times 2$ block real linear system (a saddle point system)
\begin{equation}\label{eqn:main_system}
		\underbrace{
            \begin{bmatrix} 
			 \beta I  & K + \alpha I \\
			 K + \alpha I & -\beta I
		\end{bmatrix}}_{=:\mathcal{{A}}}
		\begin{bmatrix} \mathbf{z}_1\\ \mathbf{z}_2 \end{bmatrix}=
		\begin{bmatrix} \mathbf{a}\\ \mathbf{b} \end{bmatrix},
	\end{equation} 
 where $\lambda=\alpha+\beta {\bf i}$, $\textbf{z}=\textbf{z}_1+\textbf{z}_2{\bf i}$, $\textbf{{f}}=\textbf{a}+\textbf{b}{\bf i}$, with ${\bf i}=\sqrt{-1}$ denoting the imaginary unit. Such a ``complex-to-real" equivalent transformation is not new and several preconditioners for such saddle point systems have been studied in the literature. In \cite{van2012numerical,axelsson2000real}, Schur complement based preconditioners have been proposed for the saddle point system. Recently, the preconditioned square block (PRESB) preconditioner has been proposed in \cite{axelsson2020superior,axelsson2020optimality} for the saddle point system, which is based on block triangular factorization formulas and thus easily invertible. {Nevertheless, when applying the PRESB preconditioner to \eqref{eqn:main_system}, its theoretical analysis relies on the assumption that  $K + \alpha I+\beta I$ is SPD \cite{axelsson2024efficient}. Moreover, inverting $K + \alpha I+\beta I$ to implement these preconditioners presents a challenge. When $K + \alpha I+\beta I$ is indefinite, its inversion is as difficult as to invert the original system $(K+\lambda I)\textbf{z}=\textbf{f}$. To summarize, iterative solver for \eqref{eqn:main_system} with indefinite $K + \alpha I+\beta I$ is still a challenging problem.}
 
 For any real symmetric matrix ${\bf H}$, define its absolute-value as
 \begin{equation*}
 |{\bf H}|:=Q^{\rm T}{\rm diag}(|\lambda_1|,|\lambda_2|,...,|\lambda_M|)Q,
 \end{equation*}
 where $Q^{\rm T}{\rm diag}(\lambda_1,\lambda_2,...,\lambda_M)Q$ is the unitary diagonalization of ${\bf H}$; in particular, if $\lambda_i$'s are nonnegative, we define the square root of ${\bf H}$ as 
 \begin{equation*}
 \sqrt{\bf H}:=Q^{\rm T}{\rm diag}(\sqrt{\lambda_1},\sqrt{\lambda_2},...,\sqrt{\lambda_M})Q.
 \end{equation*}
 
  Then, our absolute-value preconditioner for \eqref{eqn:main_system} is given by
 	\begin{eqnarray}\label{eqn:matrix_optipre}
 	\mathcal{\tilde{A}}:=|\mathcal{A}|=\sqrt{\mathcal{A}^T*\mathcal{A}}&=&
 		\left[ \begin{array}[c]{cc}
 			\sqrt{(K + \alpha I)^2+ \beta^2 I}&\\
 			&\sqrt{(K + \alpha I)^2+ \beta^2 I}
 		\end{array}\right].
 		\end{eqnarray}
		
	To demonstrate the effectiveness of the proposed preconditioner, we will show that the preconditioned $\mathcal{\tilde{A}}^{-1}\mathcal{{A}}$ has only two distinct eigenvalues $\pm 1$, leading to a two-step convergence of the MINRES solver (i.e., the MINRES solver for the preconditioned system will find the exact solution to \eqref{eqn:main_system} within two-iterations). In each iteration of the preconditioned MINRES solver, it requires to compute  matrix vector multiplications $\mathcal{\tilde{A}}^{-1}{\bf v}$ for some given vector ${\bf v}$. When $a({\bf x})$ is a constant, $\mathcal{\tilde{A}}^{-1}$ is diagonalizable by discrete sine transforms and such computing such matrix-vector multiplication requires $\mathcal{O}(M\log M)$ operations. When $a({\bf x})$ is non-constant, computing such matrix-vector multiplication would be expensive. To remedy this, we employ the following absolute-value based preconditioner in the case of non-constant $a({\bf x})$.
	\begin{equation}\label{eqn:preconditioner_PP}
		\mathcal{P} = \left[ \begin{array}[c]{cc}
			\sqrt{(\Bar{K} + \alpha I)^2+ \beta^2 I}&\\
			&\sqrt{(\Bar{K} + \alpha I)^2+ \beta^2 I}
		\end{array}\right],
	\end{equation}
where $\Bar{K}$ = $\gamma L$ for $\gamma=\sqrt{\check{a}\hat{a}}$. It is clear that $L$ is the discretization matrix of the constant Laplacian $-\nabla^2$, which can be diagonalized by the discrete sine transforms and its eigenvalues are larger than a positive constant $c_0$. It is important to note that when $a(x)$ is equal to 1, the preconditioner $\mathcal{P} $ is exactly the preconditioner 	$\mathcal{\tilde{A}}$ and $K=L$. {Such a preconditioner is always fast invertible no matter $K+\alpha I$ is definite or not. Moreover, we prove that the preconditioned MINRES solver with the proposed preconditioner has a convergence rate independent of matrix size under assumptions on $\alpha$ and $\beta$\footnote{See Theorem \ref{mainthm2}}, even if $K + \alpha I$ may be indefinite in such assumptions. 
Hence, the proposed preconditioning technique somewhat relaxes the restriction on definiteness of $K + \alpha I$, compared with the results in the literature. However, as mentioned above, our preconditioner relies on the uniform grid discretization. The uniform grid discretization does not exist when the domain of interested is irregular one. In the futrue, we will resort to the technique of domain decomposition to partition  the major area of an irregular domain with uniform grid so that our preconditioner is applicable to the major area.}

	The contribution of this work is twofold: first, we start with the real symmetric equivalent form of the complex-shifted Laplacian system \eqref{eqn:complex_shifted_Lap_system} and propose a novel absolute-value based preconditioner for the real symmetric system. Thanks to the real symmetric equivalent form, we can employ the short-recurrence Krylov subspace solver, MINRES, to solve the preconditioned system, which helps reducing the computational cost of iterative steps. Second, we develop a convergence theory for the preconditioned MINRES solver, showing that the  proposed solver has a matrix-size independent convergence rate.
	
	The paper is organized as follows. In Section \ref{sec:main}, we provide the convergence analyses for our proposed preconditioned solver. Numerical examples are given in Section \ref{sec:numerical} for showing the effectiveness of our proposed preconditioners. At last, the concluding remarks are given in {Section} \ref{sec:conclusions}.
	
\section{Theoretical results}\label{sec:main}
In this section, we investigative the convergence rate of MINRES solver for the preconditioned system. 
The following lemma will be exploited as framework for analyzing the convergence behavior of the proposed preconditioned MINRES method.
\begin{lemma}\label{minrescvglm}\cite[Theorem 6.13]{elman2014finite}
	Let ${\bf P}\in\mathbb{R}^{N\times N}$ and ${\bf A}\in\mathbb{R}^{N\times N}$ be a symmetric positive definite matrix and a symmetric nonsingular matrix, respectively. Suppose $\Sigma({\bf P}^{-1}{\bf A})\in[-a_1,-a_2]\cup[a_3,a_4]$ with $a_4\geq a_3>0$, $a_1\geq a_2>0$ and $a_1-a_2=a_4-a_3$. Then, the MINRES solver with ${\bf P}$ as a preconditioner for the linear system ${\bf A}{\bf x}={\bf y}\in\mathbb{R}^{N\times 1}$ has a linear convergence as follows
	\begin{equation*}
		||{\bf r}_{k}||_2\leq 2\left(\frac{\sqrt{a_1a_4}-\sqrt{a_2a_3}}{\sqrt{a_1a_4}+\sqrt{a_2a_3}}\right)^{\lfloor k/2\rfloor}||{\bf r}_{0}||_2,
	\end{equation*}
	where ${\bf r}_{k}={\bf P}^{-1}{\bf y}-{\bf P}^{-1}{\bf A}{\bf x}_k$ denotes the residual vector at the $k$th iteration with ${\bf x}_k$ ($k\geq 1$) being the $k$th iterative solution by MINRES; ${\bf x}_0$ denotes an arbitrary real-valued initial guess; $\lfloor k/2\rfloor$ denotes the integer part of $k/2$.
\end{lemma}


\subsection{Spectral analysis of the ideal preconditioner $\mathcal{\tilde{A}}^{-1}\mathcal{A}$}


The following theorem reveals the eigendecomposition structure of a class of certain $2 \times 2$ block matrices, including the matrix $\mathcal{A}$ as a special case. Note that such a block structure can be regarded as a special case in our early work on asymptotic eigenvalue distribution \cite{ferrari2019eigenvalue,hon2023preconditioned}.

\begin{theorem}\label{thm:eig_matrix_A}
	Let 
	\begin{eqnarray}\label{eqn:gen_mat_m}
		\mathcal{M} = \begin{bmatrix} 
			\theta I_{n}  & A_n^*\\
			A_n  &  -\theta I_{n}
		\end{bmatrix}
		\in \mathbb{C}^{2n \times 2n}
	\end{eqnarray} with $\theta \in \mathbb{R}$ and $A_n \in \mathbb{C}^{n \times n}$, and let $|\mathcal{M}|:=\sqrt{\mathcal{M}^2}$. Then, $|\mathcal{M}|^{-1}\mathcal{M}$ is both Hermitian and unitary. In other words, $|\mathcal{M}|^{-1}\mathcal{M}$ has only two distinct eigenvalues $\pm 1$.
\end{theorem}
\begin{proof}
	Considering the singular value decomposition of $A_n = U_n \Sigma_n V_n^*$, we have the following eigendecomposition of $\mathcal{M}$ by direct computations
	\begin{eqnarray}\label{eqn:matrix_Mat}
		\mathcal{M} &=&
		\begin{bmatrix} 
			\theta I_{n}  & V_n \Sigma_n U_n^*\\
			U_n \Sigma_n V_n^*  &  -\theta I_n
		\end{bmatrix} \\\nonumber
		&=&
		\begin{bmatrix} 
			V_n  & \\
			&  U_n
		\end{bmatrix}
		\begin{bmatrix} 
			\theta I_{n}  & \Sigma_n \\
			\Sigma_n  &  -\theta I_{n}
		\end{bmatrix}
		\begin{bmatrix} 
			V_n  & \\
			&  U_n
		\end{bmatrix}^*  \\\nonumber
		&=&\underbrace{\begin{bmatrix} 
				V_n  & \\
				&  U_n
			\end{bmatrix}
			\mathcal{\widetilde{Q}}}_{=:\mathcal{Q}}
		\begin{bmatrix} 
			\sqrt{\Sigma_n^2+\theta^2 I_{n}}  & \\
			&  -\sqrt{\Sigma_n^2+\theta^2 I_{n} }
		\end{bmatrix}
		\underbrace{
			\mathcal{\widetilde{Q}}^*
			\begin{bmatrix} 
				V_n  & \\
				&  U_n
			\end{bmatrix}^*}_{=\mathcal{Q^*}},
	\end{eqnarray} 
	where $\mathcal{\widetilde{Q}}$ is real orthogonal given by 
	\begin{eqnarray*}
		\mathcal{\widetilde{Q}}
		=
		\begin{bmatrix} 
			\Sigma_n D_1^{-1}& -\Sigma_n D_2^{-1}\\
			\big(\sqrt{\Sigma_n^2+\theta^2 I_{n}}  - \theta  I_{n}\big)D_1^{-1}   & \big(\sqrt{\Sigma_n^2+\theta^2 I_{n}}  + \theta I_{n}\big)D_2^{-1}
		\end{bmatrix}
	\end{eqnarray*}
	with
	\[
	D_1 =  \sqrt{(-\sqrt{\Sigma_n^2+\theta^2 I_{n}}  + \theta I_{n})^2 + \Sigma_n^2} 
	\]
	and
	\[
	D_2 = \sqrt{(\sqrt{\Sigma_n^2+\theta^2 I_{n}}  + \theta I_{n})^2 + \Sigma_n^2 }.
	\]
	It is obvious that both diagonal matrices $D_1$ and $D_2$ are invertible, since $\theta \neq 0$ by assumption. 
	
	With \eqref{eqn:matrix_Mat}, it is easy to see that
	\begin{eqnarray*} 
		|\mathcal{M}|&=&\sqrt{\mathcal{M}^2}\\
		&=&\sqrt{\mathcal{Q}\left[\begin{array}[c]{cc}
				\Sigma_n^2+\theta^2 I_{n}&\\
				&\Sigma_n^2+\theta^2 I_{n}
			\end{array}\right]\mathcal{Q}^{*}
		}\\
		&=&\mathcal{Q}\left[ \begin{array}[c]{cc}
			\sqrt{\Sigma_n^2+\theta^2 I_{n}}&\\
			&\sqrt{\Sigma_n^2+\theta^2 I_{n}}
		\end{array}\right]\mathcal{Q}^{*}.
	\end{eqnarray*}
	Therefore, we have
	\begin{equation*}
		|\mathcal{M}|^{-1}  \mathcal{M}=\mathcal{Q}\left[ \begin{array}[c]{cc}
			I_n&\\
			&-I_n
		\end{array}\right]\mathcal{Q}^{*}.
	\end{equation*}
	It is easy to see from the equality above that $|\mathcal{M}|^{-1}  \mathcal{M}$ is both unitary and Hermitian with $\pm 1$ as its eigenvalues.
	The proof is concluded.
\end{proof}

Let $\mathcal{A}$ and ${\tilde{\mathcal{A}}}$ be defined by \eqref{eqn:main_system} and \eqref{eqn:matrix_optipre}, respectively. Since $K + \alpha I$ is symmetric, it is straight forward to verify that $\mathcal{\tilde{A}}^{-1}\mathcal{A}$  has only two distinct eigenvalues $\pm 1$ with using Theorem \ref{thm:eig_matrix_A}. Then, with Lemma \ref{minrescvglm}, we immediately obtained the following theorem.

\begin{theorem}\label{minrescvgthm}
For any real initial guess, the MINRES solver with $\mathcal{\tilde{A}}$ as preconditioner will find the exact solution within 2 iterations.
\end{theorem}
\begin{proof}
Since the spectrum of $\mathcal{\tilde{A}}^{-1}\mathcal{A}$ is equal to $\{-1,1\}$, the convergence factor given in Lemma \ref{minrescvglm} is equal to $0$ when $k=2$. That means $||{\bf r}_2||_2=0$ with ${\bf r}_2$ being the residual vector at 2nd iteration. The result follows.
\end{proof}

\subsection{Spectral analysis of $\mathcal{P}^{-1}\mathcal{A}$}

In this subsection, we provide a theorem which supports the optimal preconditioning effect of $\mathcal{P}$.

Denote by $\rho({\bf C})$ and $\sigma({\bf C})$ respectively the spectral radius and the spectrum of a square matrix ${\bf C}$. Denote by $\lambda_{\min}({\bf H})$ and $\lambda_{\max}({\bf H})$ respectively the maximal eigenvalue and the minimal eigenvalue of a Hermitian matrix ${\bf H}$. 

For any real symmetric matrices ${\bf C}_1,{\bf C}_2\in\mathbb{R}^{k\times k}$, denote ${\bf C}_2 \succ ({\rm or} \succeq) \ {\bf C}_1$
if ${\bf C}_2-{\bf C}_1$ is positive definite (or semi-definite). Especially, we denote ${\bf C}_2 \succ ({\rm or} \succeq) \ {\bf O}$, if ${\bf C}_2$ itself is positive definite (or semi-definite).
Also, ${\bf C}_1 \prec ({\rm or} \preceq) \ {\bf C}_2$ and ${\bf O} \prec ({\rm or} \preceq) \ {\bf C}_2$  have the same meanings as those of ${\bf C}_2 \succ ({\rm or} \succeq) \ {\bf C}_1$ and ${\bf C}_2 \succ ({\rm or} \succeq) \ {\bf O}$, respectively.

\begin{lemma}\label{lemma:matsqrtoutinlm}
Let two Hermitian positive semi-definite matrices $H_1$ and $H_2$  be simultaneously diagonalizable. Then, for any non-zero vector $z$ of compatible size, it holds that
\begin{equation*}
 \frac{\sqrt{2}}{2}\leq \frac{z^{*}\sqrt{H_1^2+H_2^2}z}{z^{*}(H_1+H_2)z}\leq 1.   
\end{equation*}
\end{lemma}
\begin{proof}
 Since $H_1$ and $H_2$ are both simultaneously diagonalizable and Hermitian, we can rewrite $H_1$ and $H_2$ as
 \begin{equation*}
     H_1=Q^{*} \Lambda_1Q, \quad H_2=Q^{*} \Lambda_2Q,
 \end{equation*}
 for some unitary matrix $Q$ and diagonal matrices $\Lambda_1$ and $\Lambda_2$. Then,
 \begin{align*}
    \frac{z^{*}\sqrt{H_1^2+H_2^2}z}{z^{*}(H_1+H_2)z}&=\frac{z^{*}\sqrt{Q^{*}(\Lambda_1^2+\Lambda_2^2)Q}z}{z^{*}Q^{*}(\Lambda_1+\Lambda_2)Qz}\\
    &=\frac{z^{*}Q^{*}\sqrt{\Lambda_1^2+\Lambda_2^2}Qz}{z^{*}Q^{*}(\Lambda_1+\Lambda_2)Qz}\\
 &\stackrel{\tilde{\bf z}=Qz}{=\joinrel=\joinrel=\joinrel=}\frac{\tilde{\bf z}^{*}\sqrt{\Lambda_1^2+\Lambda_2^2}\tilde{\bf z}}{\tilde{\bf z}^{*}(\Lambda_1+\Lambda_2)\tilde{\bf z}}.
 \end{align*}
 Notice that $\Lambda_1$ and $\Lambda_2$ are both diagonal matrices with nonegative diagonal entries. For two non-negative  numbers, $c_1$ and $c_2$, it is well-known that
	\begin{align*}
		\frac{\sqrt{2}}{2}(c_1+c_2)\leq\sqrt{c_1^2+c_2^2}\leq c_1+c_2. 
	\end{align*}
 Therefore, 
 \begin{equation*}
   \frac{\sqrt{2}}{2}  \leq \frac{\tilde{\bf z}^{*}\sqrt{\Lambda_1^2+\Lambda_2^2}\tilde{\bf z}}{\tilde{\bf z}^{*}(\Lambda_1+\Lambda_2)\tilde{\bf z}}\leq 1,
 \end{equation*}
 which implies that 
 \begin{equation*}
     \frac{\sqrt{2}}{2}  \leq\frac{z^{*}\sqrt{H_1^2+H_2^2}z}{z^{*}(H_1+H_2)z}\leq 1.
 \end{equation*}
 The proof is complete.
\end{proof}
\begin{proposition}\label{prop:positivenum}
    For positive numbers $\theta_i, \eta_i(1 \leq i \leq m)$, it obviously holds that
$$
\min _{1 \leq i \leq m} \frac{\theta_i}{\eta_i} \leq\left(\sum_{i=1}^m \eta_i\right)^{-1}\left(\sum_{i=1}^m \theta_i\right) \leq \max _{1 \leq i \leq m} \frac{\theta_i}{\eta_i} .
$$
\end{proposition}

\begin{lemma}\textnormal{(See, e.g., \cite{lin2021parallel,goluvan2013})}\label{kmatproplm}
\begin{description}
\item[(i)] $\check{a}L\preceq K\preceq \hat{a}L $.
\item[(ii)] There exists a positive constant $c_0>0$ independent of $M$ such that $L\succeq c_0I$.
\end{description}	

\end{lemma}

\begin{theorem}\label{thm:preconditioner_PP}
Assume ${\check{a}} c_0 + |\beta|- |\alpha| > 0$, {and define $\gamma=\sqrt{\check{a}\hat{a}}$.}	Let $\mathcal{A},\mathcal{P} \in \mathbb{R}^{2M \times 2M}$ be defined by (\ref{eqn:main_system}) \& (\ref{eqn:preconditioner_PP}), respectively. Then, 
\begin{description}
\item[(i)] $\sigma(\mathcal{P}^{-1}\mathcal{A})\subset\left[-\mu_0,-\frac{1}{\mu_0}\right]\cup\left[\frac{1}{\mu_0},\mu_0\right]$ for $\alpha\geq 0$;
\item[(ii)]$\sigma(\mathcal{P}^{-1}\mathcal{A})\subset\left[-\tilde{\mu}_0,-\tilde{\mu}_1\right]\cup\left[\tilde{\mu}_1, \tilde{\mu}_0\right]$ for $\alpha<0$ and ${\check{a}} c_0 + |\beta|+\alpha > 0$,
\end{description}
where
\begin{align*}
\mu_0 = \sqrt{\frac{2\hat{a}}{\check{a}}},\qquad\tilde{\mu}_0=\frac{\sqrt{2}(c_0\hat{a}+|\beta|-\alpha)}{c_0\gamma+|\beta|+\alpha},\quad \tilde{\mu}_1={\frac{\sqrt{2}(c_0\check{a}+|\beta|+\alpha)}{2(c_0\gamma+|\beta|-\alpha)}}.
\end{align*} 
\end{theorem}
\begin{proof}
	By matrix similarity, we have
	\begin{equation}\label{matsimilarityspectrumeq}
		\sigma(\mathcal{P}^{-1}\mathcal{A})=\sigma(\mathcal{P}^{-\frac{1}{2}}\mathcal{A}\mathcal{P}^{-\frac{1}{2}})=\sigma(\mathcal{P}^{-\frac{1}{2}}\mathcal{\tilde{A}}^{\frac{1}{2}}\mathcal{\tilde{A}}^{-\frac{1}{2}}\mathcal{A}\mathcal{\tilde{A}}^{-\frac{1}{2}}\mathcal{\tilde{A}}^{\frac{1}{2}}\mathcal{P}^{-\frac{1}{2}}).
	\end{equation}
	Suppose $\lambda$ is an eigenvalue of $\mathcal{P}^{-\frac{1}{2}}\mathcal{\tilde{A}}^{\frac{1}{2}}\mathcal{\tilde{A}}^{-\frac{1}{2}}\mathcal{A}\mathcal{\tilde{A}}^{-\frac{1}{2}}\mathcal{\tilde{A}}^{\frac{1}{2}}\mathcal{P}^{-\frac{1}{2}}$. Then,
	\begin{equation*}
		|\lambda|\leq  ||\mathcal{P}^{-\frac{1}{2}}\mathcal{\tilde{A}}^{\frac{1}{2}}\mathcal{\tilde{A}}^{-\frac{1}{2}}\mathcal{A}\mathcal{\tilde{A}}^{-\frac{1}{2}}\mathcal{\tilde{A}}^{\frac{1}{2}}\mathcal{P}^{-\frac{1}{2}}||_2\leq ||\mathcal{P}^{-\frac{1}{2}}\mathcal{\tilde{A}}^{\frac{1}{2}}||_2||\mathcal{\tilde{A}}^{-\frac{1}{2}}\mathcal{A}\mathcal{\tilde{A}}^{-\frac{1}{2}}||_2||\mathcal{\tilde{A}}^{\frac{1}{2}}\mathcal{P}^{-\frac{1}{2}}||_2.
	\end{equation*}
	Since $\sigma(\mathcal{\tilde{A}}^{-\frac{1}{2}}\mathcal{A}\mathcal{\tilde{A}}^{-\frac{1}{2}})=\sigma(\mathcal{\tilde{A}}^{-1}\mathcal{A})\subset\{-1,1\}$ and $\mathcal{\tilde{A}}^{-\frac{1}{2}}\mathcal{A}\mathcal{\tilde{A}}^{-\frac{1}{2}}$ is Hermitian, it is clear that
	\begin{equation*}
		||\mathcal{\tilde{A}}^{-\frac{1}{2}}\mathcal{A}\mathcal{\tilde{A}}^{-\frac{1}{2}}||_2=1.
	\end{equation*}
	Then, the upper bound estimation of $|\lambda|$ reduces to
	\begin{equation*}
		|\lambda|\leq   ||\mathcal{P}^{-\frac{1}{2}}\mathcal{\tilde{A}}^{\frac{1}{2}}||_2||\mathcal{\tilde{A}}^{\frac{1}{2}}\mathcal{P}^{-\frac{1}{2}}||_2=||\mathcal{P}^{-\frac{1}{2}}\mathcal{\tilde{A}}^{\frac{1}{2}}||_2^2.
	\end{equation*}
	Notice that
	\begin{equation*}
		||\mathcal{P}^{-\frac{1}{2}}\mathcal{\tilde{A}}^{\frac{1}{2}}||_2^2=\lambda_{\max}(\mathcal{P}^{-\frac{1}{2}}\mathcal{\tilde{A}}\mathcal{P}^{-\frac{1}{2}}).
	\end{equation*}
	Let ${\bf z}$ be any eigenvector of $\mathcal{P}^{-\frac{1}{2}}\mathcal{\tilde{A}}\mathcal{P}^{-\frac{1}{2}}$ corresponding an eigenvalue $\lambda_{\max}(\mathcal{P}^{-\frac{1}{2}}\mathcal{\tilde{A}}\mathcal{P}^{-\frac{1}{2}})$. Then, it holds that
	\begin{eqnarray*}
		&&\lambda_{\max}(\mathcal{P}^{-\frac{1}{2}}\mathcal{\tilde{A}}\mathcal{P}^{-\frac{1}{2}})\\
		&=&\frac{{\bf z}^{\rm T}\mathcal{P}^{-\frac{1}{2}}\mathcal{\tilde{A}}\mathcal{P}^{-\frac{1}{2}}{\bf z}}{{\bf z}^{\rm T}{\bf z}}\\
&\stackrel{\tilde{\bf z}=\mathcal{P}^{-\frac{1}{2}}{\bf z} \neq 0}{=\joinrel=\joinrel=\joinrel=\joinrel=\joinrel=\joinrel=\joinrel=}&\frac{\tilde{\bf z}^{\rm T}\mathcal{\tilde{A}}\tilde{\bf z}}{\tilde{\bf z}^{\rm T}\mathcal{P}\tilde{\bf z}}\\	&=&\frac{\left[\tilde{\bf z}_1^{\rm T},\tilde{\bf z}_2^{\rm T}\right] {\rm blkdiag}\left(\sqrt{(K + \alpha I)^2+ \beta^2 I} ,\sqrt{(K + \alpha I)^2+ \beta^2 I} \right)\begin{bmatrix}
\tilde{\bf z}_1     \\
\tilde{\bf z}_2
\end{bmatrix}}{\left[\tilde{\bf z}_1^{\rm T},\tilde{\bf z}_2^{\rm T}\right]{\rm blkdiag}\left(\sqrt{(\Bar{K} + \alpha I)^2+ \beta^2 I} ,\sqrt{(\Bar{K} + \alpha I)^2+ \beta^2 I} \right)\begin{bmatrix}
\tilde{\bf z}_1     \\
\tilde{\bf z}_2
\end{bmatrix}}\\
&=& \frac{\tilde{\bf z}_1^{\rm T}\sqrt{(K + \alpha I)^2+ \beta^2 I}\tilde{\bf z}_1 + \tilde{\bf z}_2^{\rm T}\sqrt{(K + \alpha I)^2+ \beta^2 I}\tilde{\bf z}_2}{\tilde{\bf z}_1^{\rm T}\sqrt{(\Bar{K} + \alpha I)^2+ \beta^2 I}\tilde{\bf z}_1 + \tilde{\bf z}_2^{\rm T}\sqrt{(\Bar{K} + \alpha I)^2+ \beta^2 I}\tilde{\bf z}_2},
	\end{eqnarray*}
	which together with Proposition \ref{prop:positivenum} and Lemma \ref{lemma:matsqrtoutinlm} implies that
	\begin{align*}
		|\lambda|\leq \lambda_{\max}(\mathcal{P}^{-\frac{1}{2}}\mathcal{\tilde{A}}\mathcal{P}^{-\frac{1}{2}})&\leq  \max\limits_{\tilde{\bf z}_1\neq 0,\tilde{\bf z}_2\neq 0} \left\{ \frac{\tilde{\bf z}_1^{\rm T}\sqrt{(K + \alpha I)^2+ \beta^2 I}\tilde{\bf z}_1 }{\tilde{\bf z}_1^{\rm T}\sqrt{(\Bar{K} + \alpha I)^2+ \beta^2 I}\tilde{\bf z}_1 },  \frac{ \tilde{\bf z}_2^{\rm T}\sqrt{(K + \alpha I)^2+ \beta^2 I}\tilde{\bf z}_2}{\tilde{\bf z}_2^{\rm T}\sqrt{(\Bar{K} + \alpha I)^2+ \beta^2 I}\tilde{\bf z}_2}\right\}\\
		&=\max\limits_{\tilde{\bf z}_1\neq 0}\frac{ \tilde{\bf z}_1^{\rm T}\sqrt{(K + \alpha I)^2+ \beta^2 I}\tilde{\bf z}_1}{\tilde{\bf z}_1^{\rm T}\sqrt{(\Bar{K} + \alpha I)^2+ \beta^2 I}\tilde{\bf z}_1}\\
		&\leq \sqrt{2}\max\limits_{\tilde{\bf z}_1\neq 0}\frac{\tilde{\bf z}_1^{\rm T}(|K + \alpha I|+ |\beta| I)\tilde{\bf z}_1}{\tilde{\bf z}_1^{\rm T}(|\Bar{K} + \alpha I|+ |\beta| I)\tilde{\bf z}_1}.
	\end{align*}

By Proposition \ref{prop:positivenum} and Lemma \ref{kmatproplm}${\bf (i)}$, we have
\begin{align*}
		|\lambda| & \leq \sqrt{2} \max \left\{ \frac{\tilde{\bf z}_1^{\rm T}(K + \alpha I)\tilde{\bf z}_1}{\tilde{\bf z}_1^{\rm T}(\Bar{K} + \alpha I)\tilde{\bf z}_1}, 1\right\}\notag\\
   & \leq  \sqrt{2} \max \left\{ \frac{\tilde{\bf z}_1^{\rm T}K \tilde{\bf z}_1}{\tilde{\bf z}_1^{\rm T}\Bar{K}\tilde{\bf z}_1}, 1\right\}\notag\\
   &\leq  \sqrt{2} \max \left\{ \frac{\hat{a}\tilde{\bf z}_1^{\rm T}L \tilde{\bf z}_1}{\gamma\tilde{\bf z}_1^{\rm T}L\tilde{\bf z}_1}, 1\right\}= \sqrt{\frac{2\hat{a}}{\check{a}}}=\mu_0,\quad {\rm~if~}\alpha\geq 0.
	\end{align*}
	That implies
	\begin{equation}\label{nngalphaoutbdesti}
	\lambda\in[-\mu_0,\mu_0],\quad {\rm~if~}\alpha\geq 0.
	\end{equation}

For $\alpha<0$, Lemma \ref{kmatproplm}${\bf (ii)}$ implies that
\begin{align}\label{upperpartintermediatneq}
|\Bar{K} + \alpha I|+ |\beta| I\succeq |\Bar{K}|-|\alpha|I+|\beta|I=\bar{K}+|\beta|I-|\alpha| I&=\gamma L+|\beta|I-|\alpha| I\\
&{\succeq (\check{a}c_0+|\beta|-|\alpha|)I}.
\end{align}
Hence, 
\begin{align*}
|\lambda|&\leq \sqrt{2}\max\limits_{\tilde{\bf z}\neq 0}\frac{\tilde{\bf z}_1^{\rm T}(|K + \alpha I|+ |\beta| I)\tilde{\bf z}_1}{\tilde{\bf z}_1^{\rm T}(|\Bar{K} + \alpha I|+ |\beta| I)\tilde{\bf z}_1}\notag\\
&\leq \sqrt{2}\max\limits_{\tilde{\bf z}\neq 0}\frac{\tilde{\bf z}_1^{\rm T}(K +(|\alpha|+|\beta|) I)\tilde{\bf z}_1}{\tilde{\bf z}_1^{\rm T}(|\Bar{K} + \alpha I|+ |\beta| I)\tilde{\bf z}_1}\notag\\
&\leq \sqrt{2}\max\limits_{\tilde{\bf z}\neq 0}\frac{\tilde{\bf z}_1^{\rm T}(\hat{a}L +(|\alpha|+|\beta|) I)\tilde{\bf z}_1}{\tilde{\bf z}_1^{\rm T}(|\Bar{K} + \alpha I|+ |\beta| I)\tilde{\bf z}_1}\notag\\
&\leq  \sqrt{2}\max\limits_{\tilde{\bf z}\neq 0}\frac{\tilde{\bf z}_1^{\rm T}(\hat{a}L +(|\alpha|+|\beta|) I)\tilde{\bf z}_1}{\tilde{\bf z}_1^{\rm T}(\gamma L+|\beta|I-|\alpha| I)\tilde{\bf z}_1}\notag\\
&=\sqrt{\frac{2\hat{a}}{\check{a}}}\left[1+\max\limits_{\tilde{\bf z}\neq 0}\frac{[(\hat{a}^{-1}-\gamma^{-1})|\beta|+(\hat{a}^{-1}+\gamma^{-1})|\alpha|]\tilde{\bf z}_1^{\rm T}\tilde{\bf z}_1}{\tilde{\bf z}_1^{\rm T}[L+\gamma^{-1}(|\beta|-|\alpha|) I]\tilde{\bf z}_1}\right]\notag\\
&\leq \sqrt{\frac{2\hat{a}}{\check{a}}}\left[1+\max\limits_{\tilde{\bf z}\neq 0}\frac{[(\hat{a}^{-1}-\gamma^{-1})|\beta|+(\hat{a}^{-1}+\gamma^{-1})|\alpha|]\tilde{\bf z}_1^{\rm T}\tilde{\bf z}_1}{\tilde{\bf z}_1^{\rm T}[c_0I+\gamma^{-1}(|\beta|-|\alpha|) I]\tilde{\bf z}_1}\right]\notag\\
&=\frac{\sqrt{2}(|\beta|+|\alpha|+c_0\hat{a})}{\gamma c_0+|\beta|-|\alpha|}\\
&{=\frac{\sqrt{2}(|\beta|-\alpha+c_0\hat{a})}{\gamma c_0+|\beta|+\alpha}}
=\tilde{\mu}_0,\quad {\rm~if~}\alpha<0{\rm~and~} {\gamma} c_0+|\beta|-|\alpha|>0,
\end{align*}
where the 3rd inequality is obtained from Lemma \ref{kmatproplm}${\bf (i)}$; the 4th inequality is obtained from \eqref{upperpartintermediatneq}; the 5th inequality is obtained from Lemma \ref{kmatproplm}${\bf (ii)}$.
	That implies
\begin{equation}\label{ngalphaoutbdesti}
	\lambda\in[-\tilde{\mu}_0,\tilde{\mu}_0],\quad {\rm~if~}\alpha<0{\rm~and~} {\gamma} c_0+|\beta|-|\alpha|>0.
\end{equation}

 It remains to estimate the lower bound of $|\lambda|$.
 Denote by $\sigma_{\min}({\bf C})$, the minimal singular value of a matrix ${\bf C}$. As $\mathcal{P}^{-\frac{1}{2}}\mathcal{\tilde{A}}^{\frac{1}{2}}\mathcal{\tilde{A}}^{-\frac{1}{2}}\mathcal{A}\mathcal{\tilde{A}}^{-\frac{1}{2}}\mathcal{\tilde{A}}^{\frac{1}{2}}\mathcal{P}^{-\frac{1}{2}}$ is Hermitian, it is clear that
	\begin{align*}
		|\lambda|&\geq \sigma_{\min}(\mathcal{P}^{-\frac{1}{2}}\mathcal{\tilde{A}}^{\frac{1}{2}}\mathcal{\tilde{A}}^{-\frac{1}{2}}\mathcal{A}\mathcal{\tilde{A}}^{-\frac{1}{2}}\mathcal{\tilde{A}}^{\frac{1}{2}}\mathcal{P}^{-\frac{1}{2}})\\
		&=\frac{1}{||(\mathcal{P}^{-\frac{1}{2}}\mathcal{\tilde{A}}^{\frac{1}{2}}\mathcal{\tilde{A}}^{-\frac{1}{2}}\mathcal{A}\mathcal{\tilde{A}}^{-\frac{1}{2}}\mathcal{\tilde{A}}^{\frac{1}{2}}\mathcal{P}^{-\frac{1}{2}})^{-1}||_2}\\
		&=\frac{1}{||\mathcal{P}^{\frac{1}{2}}\mathcal{\tilde{A}}^{-\frac{1}{2}}\mathcal{\tilde{A}}^{\frac{1}{2}}\mathcal{A}^{-1}\mathcal{\tilde{A}}^{\frac{1}{2}}\mathcal{\tilde{A}}^{-\frac{1}{2}}\mathcal{P}^{\frac{1}{2}}||_2}\\
		&\geq \frac{1}{||\mathcal{P}^{\frac{1}{2}}\mathcal{\tilde{A}}^{-\frac{1}{2}}||_2||\mathcal{\tilde{A}}^{\frac{1}{2}}\mathcal{A}^{-1}\mathcal{\tilde{A}}^{\frac{1}{2}}||_2||\mathcal{\tilde{A}}^{-\frac{1}{2}}\mathcal{P}^{\frac{1}{2}}||_2}.
	\end{align*}
	As $\mathcal{\tilde{A}}^{\frac{1}{2}}\mathcal{A}^{-1}\mathcal{\tilde{A}}^{\frac{1}{2}}$ is unitary, $||\mathcal{\tilde{A}}^{\frac{1}{2}}\mathcal{A}^{-1}\mathcal{\tilde{A}}^{\frac{1}{2}}||_2=1$. Then,
	\begin{equation*}
		|\lambda|\geq  \frac{1}{||\mathcal{P}^{\frac{1}{2}}\mathcal{\tilde{A}}^{-\frac{1}{2}}||_2||\mathcal{\tilde{A}}^{-\frac{1}{2}}\mathcal{P}^{\frac{1}{2}}||_2}=\frac{1}{||\mathcal{P}^{\frac{1}{2}}\mathcal{\tilde{A}}^{-\frac{1}{2}}||_2^2}=\frac{1}{\lambda_{\max}(\mathcal{\tilde{A}}^{-\frac{1}{2}}\mathcal{P}\mathcal{\tilde{A}}^{-\frac{1}{2}})}.
	\end{equation*}
	Moreover,
	\begin{align*}
		&\lambda_{\max}(\mathcal{\tilde{A}}^{-\frac{1}{2}}\mathcal{P}\mathcal{\tilde{A}}^{-\frac{1}{2}})\\
		&=\max\limits_{{\bf w}\neq 0}\frac{{\bf w}^{\rm T}\mathcal{\tilde{A}}^{-\frac{1}{2}}\mathcal{P}\mathcal{\tilde{A}}^{-\frac{1}{2}}{\bf w}}{{\bf w}^{\rm T}{\bf w}}\\
		&\stackrel{\tilde{\bf w}=\mathcal{\tilde{A}}^{-\frac{1}{2}}{\bf w}\neq 0}{=\joinrel=\joinrel=\joinrel=\joinrel=\joinrel=\joinrel=\joinrel=}\max\limits_{\tilde{\bf w}\neq 0}\frac{\tilde{\bf w}^{\rm T}\mathcal{P}\tilde{\bf w}}{\tilde{\bf w}^{\rm T}\mathcal{\tilde{A}}\tilde{\bf w}}\\	
		&\stackrel{\tilde{\bf w}^{\rm T}=[\tilde{\bf w}_1^{\rm T},\tilde{\bf w}_2^{\rm T}], \tilde{\bf w}_1,\tilde{\bf w}_2\in\mathbb{R}^{M\times 1}}{=\joinrel=\joinrel=\joinrel=\joinrel=\joinrel=\joinrel=\joinrel=\joinrel=\joinrel=\joinrel=\joinrel=\joinrel=\joinrel=\joinrel=}\\
		&=\max\limits_{\tilde{\bf w}_1,\tilde{\bf w}_2\neq 0}\frac{[\tilde{\bf w}_1^{\rm T},\tilde{\bf w}_2^{\rm T}]{\rm blkdiag}\left(\sqrt{(\Bar{K} + \alpha I)^2+ \beta^2 I},\sqrt{(\Bar{K} + \alpha I)^2+ \beta^2 I} \right)\begin{bmatrix}
\tilde{\bf w}_1     \\
\tilde{\bf w}_2
\end{bmatrix}}{[\tilde{\bf w}_1^{\rm T},\tilde{\bf w}_2^{\rm T}]{\rm blkdiag}\left(\sqrt{(K + \alpha I)^2+ \beta^2 I},\sqrt{(K + \alpha I)^2+ \beta^2 I} \right)\begin{bmatrix}
\tilde{\bf w}_1     \\
\tilde{\bf w}_2
\end{bmatrix}}\\
&= \max\limits_{\tilde{\bf w}_1,\tilde{\bf w}_2\neq 0}\frac{\tilde{\bf w}_1^{\rm T}\sqrt{(\Bar{K} + \alpha I)^2+ \beta^2 I}\tilde{\bf w}_1 + \tilde{\bf w}_2^{\rm T}\sqrt{(\Bar{K} + \alpha I)^2+ \beta^2 I}\tilde{\bf w}_2}{\tilde{\bf w}_1^{\rm T}\sqrt{(K + \alpha I)^2+ \beta^2 I}\tilde{\bf w}_1 + \tilde{\bf z}_2^{\rm T}\sqrt{(K + \alpha I)^2+ \beta^2 I}\tilde{\bf w}_2}\\
&\leq\max\limits_{\tilde{\bf w}_1,\tilde{\bf w}_2\neq 0}\max \left\{ \frac{\tilde{\bf w}_1^{\rm T}\sqrt{(\Bar{K} + \alpha I)^2+ \beta^2 I}\tilde{\bf w}_1 }{\tilde{\bf w}_1^{\rm T}\sqrt{(K + \alpha I)^2+ \beta^2 I}\tilde{\bf w}_1 },  \frac{ \tilde{\bf w}_2^{\rm T}\sqrt{(\Bar{K} + \alpha I)^2+ \beta^2 I}\tilde{\bf w}_2}{\tilde{\bf w}_2^{\rm T}\sqrt{(K + \alpha I)^2+ \beta^2 I}\tilde{\bf w}_2}\right\}\\
&=\max\limits_{\tilde{\bf w}_1\neq 0}\frac{\tilde{\bf w}_1^{\rm T}\sqrt{(\Bar{K} + \alpha I)^2+ \beta^2 I}\tilde{\bf w}_1 }{\tilde{\bf w}_1^{\rm T}\sqrt{(K + \alpha I)^2+ \beta^2 I}\tilde{\bf w}_1 }\leq \sqrt{2}\max\limits_{\tilde{\bf w}_1\neq 0}\frac{\tilde{\bf w}_1^{\rm T}(|\Bar{K} + \alpha I|+ |\beta| I)\tilde{\bf w}_1 }{\tilde{\bf w}_1^{\rm T}(|K + \alpha I|+ |\beta| I)\tilde{\bf w}_1 },
	\end{align*}
	where the last inequality comes from Lemma \ref{lemma:matsqrtoutinlm}. That means
	\begin{equation*}
|\lambda|\geq \frac{1}{\sqrt{2}\max\limits_{\tilde{\bf w}_1\neq 0}\frac{\tilde{\bf w}_1^{\rm T}(|\Bar{K} + \alpha I|+ |\beta| I)\tilde{\bf w}_1 }{\tilde{\bf w}_1^{\rm T}(|K + \alpha I|+ |\beta| I)\tilde{\bf w}_1 }}=\frac{\sqrt{2}}{2\max\limits_{\tilde{\bf w}_1\neq 0}\frac{\tilde{\bf w}_1^{\rm T}(|\Bar{K} + \alpha I|+ |\beta| I)\tilde{\bf w}_1 }{\tilde{\bf w}_1^{\rm T}(|K + \alpha I|+ |\beta| I)\tilde{\bf w}_1 }}.
	\end{equation*}
By Proposition \ref{prop:positivenum} and Lemma \ref{kmatproplm}, we have
\begin{align*}
|\lambda|&\geq \frac{\sqrt{2}}{2\max\left\{\max\limits_{\tilde{\bf w}_1\neq 0}\frac{\tilde{\bf w}_1^{\rm T}\Bar{K} \tilde{\bf w}_1 }{\tilde{\bf w}_1^{\rm T}K  \tilde{\bf w}_1 },1\right\}}\\
&\geq \frac{\sqrt{2}}{2\max\left\{\max\limits_{\tilde{\bf w}_1\neq 0}\frac{\tilde{\bf w}_1^{\rm T}\Bar{K} \tilde{\bf w}_1 }{\tilde{\bf w}_1^{\rm T}\check{a}L  \tilde{\bf w}_1 },1\right\}}=\frac{1}{2}\sqrt{\frac{2\check{a}}{\hat{a}}}=\frac{1}{\mu_0},\quad {\rm ~if~}\alpha\geq 0,
\end{align*}
which together with \eqref{nngalphaoutbdesti} implies that
\begin{equation*}
\lambda\in\left[-\mu_0,-\frac{1}{\mu_0}\right]\cup\left[\frac{1}{\mu_0},\mu_0\right],\quad {\rm~if~}\alpha\geq 0.
\end{equation*}

For $\alpha<0$, Lemma \ref{kmatproplm}${\bf (ii)}$ implies that
\begin{align}\label{lowerpartintermediatneq}
	|K + \alpha I|+ |\beta| I\succeq |K|-|\alpha|I+|\beta|I\succeq\check{a}L+|\beta|I-|\alpha| I\succeq (\check{a}c_0+|\beta|-|\alpha|)I.
\end{align}
Hence,
\begin{align*}
|\lambda|&\geq \frac{\sqrt{2}}{2\max\limits_{\tilde{\bf w}_1\neq 0}\frac{\tilde{\bf w}_1^{\rm T}(|\Bar{K} + \alpha I|+ |\beta| I)\tilde{\bf w}_1 }{\tilde{\bf w}_1^{\rm T}(|K + \alpha I|+ |\beta| I)\tilde{\bf w}_1 }}\\
&\geq \frac{\sqrt{2}}{2\max\limits_{\tilde{\bf w}_1\neq 0}\frac{\tilde{\bf w}_1^{\rm T}(|\Bar{K} + \alpha I|+ |\beta| I)\tilde{\bf w}_1 }{\tilde{\bf w}_1^{\rm T}(\check{a}L+|\beta|I-|\alpha| I)\tilde{\bf w}_1 }}\\
&\geq \frac{\sqrt{2}}{2\max\limits_{\tilde{\bf w}_1\neq 0}\frac{\tilde{\bf w}_1^{\rm T}(\Bar{K} + |\alpha |I+ |\beta| I)\tilde{\bf w}_1 }{\tilde{\bf w}_1^{\rm T}[\check{a}L+|\beta|I-|\alpha| I]\tilde{\bf w}_1 }}\\
&=\frac{\sqrt{2\check{a}}}{2\sqrt{\hat{a}}\left[1+\max\limits_{\tilde{\bf w}_1\neq 0}\frac{[(\gamma^{-1}+\check{a}^{-1})|\alpha|+(\gamma^{-1}-\check{a}^{-1})|\beta|]\tilde{\bf w}_1^{\rm T}\tilde{\bf w}_1 }{\tilde{\bf w}_1^{\rm T}[L+\check{a}^{-1}(|\beta|-|\alpha|)I]\tilde{\bf w}_1 }\right]}\\
&\geq \frac{\sqrt{2\check{a}}}{2\sqrt{\hat{a}}\left[1+\max\limits_{\tilde{\bf w}_1\neq 0}\frac{[(\gamma^{-1}+\check{a}^{-1})|\alpha|+(\gamma^{-1}-\check{a}^{-1})|\beta|]\tilde{\bf w}_1^{\rm T}\tilde{\bf w}_1 }{\tilde{\bf w}_1^{\rm T}[c_0I+\check{a}^{-1}(|\beta|-|\alpha|)I]\tilde{\bf w}_1 }\right]}\\
&={\frac{\sqrt{2}(c_0\check{a}+|\beta|-|\alpha|)}{2(c_0\gamma+|\beta|+|\alpha|)} }\\
&={\frac{\sqrt{2}(c_0\check{a}+|\beta|+\alpha)}{2(c_0\gamma+|\beta|-\alpha)} }=\tilde{\mu}_1,\quad {\rm~if~}\alpha<0{\rm~and~}{c_0\check{a}+|\beta|-|\alpha|>0},
\end{align*}
where the second inequality comes from \eqref{lowerpartintermediatneq}; the 4th inequality comes from Lemma \ref{kmatproplm}${\bf (ii)}$. The inequalities above together with \eqref{ngalphaoutbdesti} indicates that
\begin{equation*}
	\lambda\in\left[-\tilde{\mu}_0,-\tilde{\mu}_1\right]\cup\left[\tilde{\mu}_1,\tilde{\mu}_0\right],\quad {\rm~if~}\alpha<0{\rm~and~}c_0\check{a}+|\beta|+\alpha>0.
\end{equation*}
The proof is complete.
\end{proof}

With Lemma \ref{minrescvglm} and Theorem \ref{thm:preconditioner_PP}, we immediately obtain the following theorem
\begin{theorem}\label{mainthm2}
Let ${\bf r}_k$ denote the residual vector at $k$th preconditioned MINRES iteration with $\mathcal{P}$ as preconditioner for solving the system \eqref{eqn:main_system} for $k\geq 1$. Let
\begin{equation*}
{\bf r}_0=\mathcal{P}^{-1}\left[\begin{array}[c]{c}
{\bf a}\\
{\bf b}
\end{array}\right]-\mathcal{A}{\bf z}_0,
\end{equation*} 
denote the intial residual vector with arbitrary real initial guess ${\bf z}_0$. Then,
\begin{description}
\item[(i)] \begin{equation*}
||{\bf r}_k||_2\leq 2\theta_1^{\lfloor k/2\rfloor}||{\bf r}_0||_2,\quad {\rm~if~}\alpha\geq 0,
\end{equation*}
where the convergence factor $\theta_1:=\frac{\mu_0^2-1}{\mu_0^2+1}\in(0,1)$ with $\mu_0$ defined in Theorem \ref{thm:preconditioner_PP} is independent of $M$.
\item[(ii)]  \begin{equation*}
||{\bf r}_k||_2\leq 2\theta_2^{\lfloor k/2\rfloor}||{\bf r}_0||_2,\quad {\rm~if~}\alpha<0{\rm~and~}c_0\check{a}+|\beta|+\alpha>0,
\end{equation*}
where the convergence factor $\theta_2:=\frac{\tilde{\mu}_0-\tilde{\mu}_1}{\tilde{\mu}_0+\tilde{\mu}_1}\in(0,1)$ with $\tilde{\mu}_0$ and $\tilde{\mu}_1$ defined in Theorem \ref{thm:preconditioner_PP} is independent of $M$.
\end{description}
\end{theorem}

\subsection{Implementations}\label{implementation}

We first discuss the computation of the matrix-vector product $\mathcal{A}\mathbf{v}$ for any given vector $\mathbf{v}$. Since $\mathcal{A}$ is a sparse matrix (due to $K \in \mathbf{C}^{M \times M}$ being sparse), computing $\mathcal{A}\mathbf{v}$ only requires linear complexity of $\mathcal{O}(M)$. 

For inversion of the preconditioner, we adopt different strategies in different cases.

{\bf {Constant Laplacian $-\nabla^2$ on rectangular domain and uniform grid:}} In such case, we employ finite difference method or finite element method on uniform grid for discretization of {Laplace operator. As a result, $K$ is diagonalizable by the discrete sine transform, i.e., $K = W \Lambda W^T$, where $W$ denotes a real symmetric orthogonal discrete sine transform;  $\Lambda$ is diagonal matrix with positive diagonal elements. Since $K$ is simultaneously diagonalizable by $W$, it is clear that $\tilde{\mathcal{A}}$ is fast diagonalizable.} Especially,
\begin{eqnarray*}
\tilde{\mathcal{A}}&=
&\left[ \begin{array}[c]{cc}
	W&\\
	&W
\end{array}\right]\\
&& \times
\left[ \begin{array}[c]{cc}
	\sqrt{(\Lambda+\alpha I)^2+\beta^2 I}&\\
	&\sqrt{(\Lambda+\alpha I)^2+\beta^2 I}
\end{array}\right]\\
&& \times
\left[ \begin{array}[c]{cc}
	W&\\
	&W
\end{array}\right].
\end{eqnarray*}
Hence, computing $\tilde{\mathcal{A}}^{-1}{\bf v}$ for a given vector ${\bf v}$ requires only $\mathcal{O}(M\log M)$ operations in such case. {As $\tilde{\mathcal{A}}$ is already fast invertible in such case. Therefore, in such case, we adopt the ideal preconditioner $\tilde{\mathcal{A}}$ for solving the system \eqref{eqn:main_system}.}

{\bf {Variable Laplacian $-\nabla(a({{\bf x}})\nabla)$ on rectangular domain and uniform grid:}} 
Our preconditioning technique begins with replacing $a(\mathbf{x})$ by a constant $\gamma$ ($\gamma$ = $\sqrt{\check{a}\hat{a}}$) to obtain the constant-coefficient matrix $\Bar{K}$, where $\Bar{K}$ = $\gamma L$. Note that $L$ is the discretization matrix of the constant Laplacian $-\nabla^2$. In such case, we adopt $\mathcal{P}$ as preconditioner. Similarly, computing $\mathcal{P}^{-1}{\bf v}$ for a given vector ${\bf v}$ requires only $\mathcal{O}(M\log M)$ operations.

\section{Numerical experiments}\label{sec:numerical}
{
To demonstrate the efficiency of the proposed preconditioning technique, we in this section test and compare its performance with several state-of-the-art solvers, including the multigrid methods with Jacobi smoother and Vanka-type smoother \cite{he2022vanka} , PRESB preconditioning technique \cite{axelsson2014comparison} and PMHSS preconditioning technique \cite{bai2007successive}. As both the PRESB and PMHSS preconditioners are non-Hermitian preconditioner, we apply GMRES solver to solve the preconditioned systems by the two preconditioners. The implementation of  the PRESB preconditioner and PMHSS preconditioner involves solving linear systems with $K + (\alpha+\beta) I$ as coefficient matrix, which is challenging when  $K + (\alpha+\beta) I$ is indefinite.} {To address this issue, we apply the ILU (Incomplete LU Decomposition) method to invert this part.}

{The experiments were conducted using GNU Octave 8.2.0 on a Dell R640 server with dual Xeon Gold 6246R 16-Cores 3.4 GHz CPUs and 512GB RAM running Ubuntu $20.04$ LTS. The CPU time in seconds was measured using the built-in \textbf{tic/toc} functions in Octave. The MINRES solver and GMRES solver were implemented using the built-in \textbf{minres} function and \textbf{gmres} function, respectively in Octave. {For the PRESB and PMHSS methods, $K + (\alpha+\beta) I$ is of incomplete factorization type and is obtained
via Matlab’s function \textbf{ilu}.} The exact solution was generated using the built-in \textbf{randn} function in Octave, where the exact solution is represented as $randn(M) + randn(M)*$ \textbf{i}. Here, randn denotes a random number drawn from a standard normal distribution, and \textbf{i} represents the imaginary unit. Additionally, for the MINRES and GMRES solver, we employed a zero initial guess and set a stopping tolerance of $10^{-8}$ based on the reduction in relative residual norms, unless otherwise specified. For the remaining parameter settings, we followed the same configuration as described in \cite{he2022vanka}.} {When employing the PRESB and PMHSS preconditioners, we utilized GMRES(50) method as the iterative scheme. The iteration was terminated when the residual norm was reduced by a factor of $10^{-8}$. }

In the tables given in the following examples, ``Iter" stands for the iteration number by required by MINRES for convergence and ``CPU'' indicates the CPU time estimated in seconds for convergence using the Octave built-in function \textbf{tic/toc}.

\begin{example}\label{ex:con1} {\bf Constant Laplacian case: a(x)=1.}
\end{example}
In the first example, {we tested the saddle point system \eqref{eqn:main_system} using MINRES with the proposed preconditioner $\mathcal{\tilde{A}}$ in \eqref{eqn:matrix_optipre} and the constant-coefficient complex-shifted Laplacian system \eqref{eqn:complex_shifted_Lap_system} by the MGM method with the Vanka-type smoother proposed in \cite{he2022vanka}, respectively.} Table \ref{table_2D_exp1_con_A_iter} (a) shows the iteration and CPU time of MINRES when the preconditioner $\mathcal{\tilde{A}}$ is applied and Table \ref{table_2D_exp1_con_A_iter} (b) shows the numerical results of Multigrid with Vanka-type solver. As we know, $\lambda=\alpha+\beta$ \textbf{i}. From Table \ref{table_2D_exp1_con_A_iter}, it can be observed that regardless of whether the complex shift $\lambda$ is positive or negative, there is no significant variation in the numerical results.

\begin{table}[htbp]\Huge\label{ex:con}
\begin{subtable}[t]{\textwidth}
		\caption{MINRES with the optimal preconditioner $\tilde{\mathcal{A}}$}
\centering
\resizebox{\textwidth}{!}{
		\begin{tabular}{|cc|cc|cc|cc|cc|cc|cc|}
			\hline
			\multicolumn{2}{|c|}{$\textrm{tol} = 10^{-8}$} & \multicolumn{2}{c|}{$(\alpha, \beta)=(100,100)$} & \multicolumn{2}{c|}{$(\alpha, \beta)=(-100,-100)$} & \multicolumn{2}{c|}{$(\alpha, \beta)=(100,-100)$} & \multicolumn{2}{c|}{$(\alpha, \beta)=(-100,100)$}&\multicolumn{2}{c|}{$(\alpha, \beta)=(-100,1)$}&\multicolumn{2}{c|}{$(\alpha, \beta)=(1,-100)$} \\ \hline
			${M}$                                 & DoF      & Iter                   & CPU    & Iter                      & CPU       & Iter                      & CPU       & Iter                      & CPU
   & Iter                      & CPU
   & Iter                      & CPU\\ \hline
			\multicolumn{1}{|c|}{$2^{6} - 1$}   & 7938   &    2                   &  0.010  &    2                      &   0.0072    &   2                       &  0.0074     &    2                      &    0.0073 & 2&0.011 & 2&0.016   \\ \hline
			\multicolumn{1}{|c|}{$2^{8} - 1$}  & 130050  &    2                   & 0.054   &    2                      &   0.060    &     2                     &   0.057    &     2                     &   0.054& 2&0.071 & 2& 0.060   \\ \hline
			\multicolumn{1}{|c|}{$2^{10} - 1$}  &  2093058 &   2                    & 0.89   &  2                        &  0.82     &    2                      &   0.83    &    2                      &  0.83 & 2&0.99 & 2&0.93   \\ \hline
			\multicolumn{1}{|c|}{$2^{12} - 1$}  & 33538050 &  2                     &  20.65 &   2                       &  22.71    &     2                     &  20.71    &    2                      &  20.55 & 2&23.22 & 2&23.09 \\ \hline
			\multicolumn{1}{|c|}{$2^{14} - 1$}  & 536805378 &   2                    &  322.73  &     2                     &  335.83    &     2                     &  356.19    &    2                      &  334.66 & 2&353.76 & 2&355.89  \\ \hline
		\end{tabular}
  }
	
		\end{subtable}

\begin{subtable}[t]{\textwidth}
		\caption{Multigrid with Vanka-type solver}
\centering
\resizebox{\textwidth}{!}{
		\begin{tabular}{|cc|cc|cc|cc|cc|cc|cc|}
			\hline
			\multicolumn{2}{|c|}{$\textrm{tol} = 10^{-8}$} & \multicolumn{2}{c|}{$(\alpha, \beta)=(100,100)$} & \multicolumn{2}{c|}{$(\alpha, \beta)=(-100,-100)$} & \multicolumn{2}{c|}{$(\alpha, \beta)=(100,-100)$} & \multicolumn{2}{c|}{$(\alpha, \beta)=(-100,100)$} & \multicolumn{2}{c|}{$(\alpha, \beta)=(-100,1)$}& \multicolumn{2}{c|}{$(\alpha, \beta)=(1,-100)$}\\ \hline
			${M}$                                 & DoF      & Iter                   & CPU    & Iter                      & CPU       & Iter                      & CPU       & Iter                      & CPU  & Iter                      & CPU& Iter                      & CPU     \\ \hline
			\multicolumn{1}{|c|}{$2^{6} - 1$}   & 3969   &    14                   &  0.052    &   14                       &  0.051     &                    14     &   0.051    &     14                     &   0.051 & 14&0.078 & 14&0.054   \\ \hline
			\multicolumn{1}{|c|}{$2^{8} - 1$}  & 65025  &  14                     & 0.31   &   14                       &  0.31     &  14                        &  0.31      &    14                      &   0.31 & 14&0.31 & 14&0.29  \\ \hline
			\multicolumn{1}{|c|}{$2^{10} - 1$}  & 1046529  &  14                     &  3.80   &  14 & 3.69          &  14                       &  3.74   &    14                      &  3.73 & 14&3.80 & 14&3.70  \\ \hline
			\multicolumn{1}{|c|}{$2^{12} - 1$}  & 16769025 & 14                      & 69.16  & 14  & 69.25           &  14                       & 69.02   &    14                     &   69.57 & 14&70.01 & 14&69.41 \\ \hline
			\multicolumn{1}{|c|}{$2^{14} - 1$}  & 268402689 &   14                     & 1172.2  &   14         & 1173.3   &              14            &  1157.0    &     14                     &  1167.9 & 14&1183.3 & 14&1145.8  \\ \hline
		\end{tabular}
  }
	
\end{subtable}
	\caption{Numbers of iterations and CPU times for Example \ref{ex:con1}}
\label{table_2D_exp1_con_A_iter}
\end{table}


Observing the data for Table \ref{table_2D_exp1_con_A_iter} for a wide range of degrees of freedom (DoF), it can be concluded that the optimal preconditioner in MINRES achieves the anticipated two-step convergence. In contrast, the Vanka-type solver in Multigrid demonstrates higher iteration numbers and CPU times compared to our proposed preconditioner. Specifically, the iteration number of the Vanka-type solver is approximately seven times that of our proposed preconditioner, while the CPU time is approximately three to five times longer than our method for the same degree of freedom (DoF). Regardless of the selection of $(\alpha, \beta)$, the stability of the iteration count in MINRES with the preconditioner $\tilde{\mathcal{A}}$ implies that its convergence rate remains optimal even as the grids are refined. This observation provides support for the theoretical results presented in Theorem \ref{thm:eig_matrix_A} and demonstrates the robustness of the proposed solver.
In summary, the results from the provided data demonstrate that the MINRES method with the optimal preconditioner is superior to the Multigrid method with a Vanka-type solver in terms of both iteration number and CPU time for the Constant Laplacian case ((a(x)=1)).

\begin{example}\label{ex:var} {\bf Variable Laplacian case: a(x)=$(20+ x_1^2)(20+ x_2^2)$.}
\end{example}
In the second example, as previously mentioned, {we performed experiments to demonstrate the versatility of our approach. Specifically, we conducted tests on the variable-coefficient complex-shifted Laplacian system \eqref{eqn:complex_shifted_Lap_system} using the MGM method with the Jacobi smoother and the GMRES method with the PMHSS preconditioner. For the saddle point system \eqref{eqn:main_system}, we evaluated the performance of the MINRES method with the proposed preconditioner $\mathcal{P}$ given in \eqref{eqn:preconditioner_PP}, as well as the GMRES method with the PRESB preconditioner.} The iteration count and CPU time of MINRES with the preconditioner $\mathcal{P}$ are presented in Table \ref{table_2D_exp1_var_A_iter} (a), while Table \ref{table_2D_exp1_var_A_iter} (b) displays the numerical results obtained with the Multigrid method using the Jacobi solver. {The numerical results shown in the remaining Tables \ref{table_2D_exp1_var_A_iter} (c) and \ref{table_2D_exp1_var_A_iter} (d) are from GMRES(50) with the PRESB preconditioner and GMRES(50) with the PMHSS preconditioner, respectively.} Similarly, across a wide range of $\alpha$ and $\beta$ values, {we observed that the proposed preconditioner $\mathcal{P}$ consistently outperformed the Jacobi solver and the remaining two GMRES solvers in terms of CPU time and iteration counts.} 

From Table \ref{table_2D_exp1_var_A_iter} (a), it can be observed that the iteration count of our proposed method remains stable at 14, even when the matrix size is twice that of the Multigrid method. This demonstrates our advantage in terms of computational efficiency. In contrast, the Jacobi solver in Multigrid exhibits significantly higher iteration numbers and CPU times for the same spatial grid ${M}$. Specifically, the iteration count of the Jacobi solver is more than twice that of our proposed preconditioner, and the CPU time greatly exceeds that of our method. Constructing a sparse variable-coefficient Laplacian matrix also requires a considerable amount of time. { For Tables \ref{table_2D_exp1_var_A_iter} (c) and (d), although the numerical results of GMRES with the PMHSS preconditioner are significantly better than those of GMRES with the PRESB preconditioner, neither approach matches the efficacy of the preconditioner we have proposed.}

    It is worth noting that when ${M} \geq 2^{10}-1$, the runtime of program exceeds 3600 seconds, and we chose to terminate the program at that point.
\begin{table}[h!]\Huge\label{ex:con}
\begin{subtable}[t]{\textwidth}{
		\caption{MINRES with the optimal preconditioner $\mathcal{P}$}
\resizebox{\textwidth}{!}{
		\begin{tabular}{|cc|cc|cc|cc|cc|cc|cc|}
			\hline
			\multicolumn{2}{|c|}{$\textrm{tol} = 10^{-8}$} & \multicolumn{2}{c|}{$(\alpha, \beta)=(-600,150)$} & \multicolumn{2}{c|}{$(\alpha, \beta)=(-100,-25)$} & \multicolumn{2}{c|}{$(\alpha, \beta)=(100,-100)$} & \multicolumn{2}{c|}{$(\alpha, \beta)=(-100,100)$}& \multicolumn{2}{c|}{$(\alpha, \beta)=(-100,1)$} & \multicolumn{2}{c|}{$(\alpha, \beta)=(1,-100)$} \\ \hline
			$M$                                 & DoF      & Iter                   & CPU    & Iter                      & CPU       & Iter                      & CPU       & Iter                      & CPU & Iter                      & CPU       & Iter                      & CPU       \\ \hline
			\multicolumn{1}{|c|}{$2^{6} - 1$}   & 7938   &    14                   &  0.026  &    14                      &   0.026    &   14                       &  0.025     &    14                      &    0.028 & 14& 0.030 & 14& 0.029     \\ \hline
			\multicolumn{1}{|c|}{$2^{8} - 1$}  & 130050  &    14                   & 0.39   &    14                      &   0.36    &     14                     &   0.24    &     14                     &   0.25 & 14&0.33 & 14&0.29    \\ \hline
			\multicolumn{1}{|c|}{$2^{10} - 1$}  &  2093058 &   14                    & 4.59   &  14                        &  4.67      &    14                      &   4.14    &    14                      &  4.47 & 14&4.13 & 14&4.34    \\ \hline
			\multicolumn{1}{|c|}{$2^{12} - 1$}  & 33538050 &  14                     &  114.74 &   14                       &  114.66    &     14                     &  110.62    &    14                      &  109.70 & 14&111.25 & 14&110.60   \\ \hline
			\multicolumn{1}{|c|}{$2^{14} - 1$}  & 536805378 &   14                    &  1670.2  &     14                     &  1666.5    &     14                     &  1680.3    &    14                     &  1683.6  & 14&1715.0 & 14&1777.1  \\ \hline
		\end{tabular}
  }}
	
		\end{subtable}

\begin{subtable}[t]{\textwidth}{
		\caption{Multigrid with the Jacobi smoother}
\resizebox{\textwidth}{!}{
		\begin{tabular}{|cc|cc|cc|cc|cc|cc|cc|}
			\hline
			\multicolumn{2}{|c|}{$\textrm{tol} = 10^{-8}$} & \multicolumn{2}{c|}{$(\alpha, \beta)=(-600,150)$} & \multicolumn{2}{c|}{$(\alpha, \beta)=(-100,-25)$} & \multicolumn{2}{c|}{$(\alpha, \beta)=(100,-100)$} & \multicolumn{2}{c|}{$(\alpha, \beta)=(-100,100)$}& \multicolumn{2}{c|}{$(\alpha, \beta)=(-100,1)$} & \multicolumn{2}{c|}{$(\alpha, \beta)=(1,-100)$} \\ \hline
			$M$                                 & DoF      & Iter                   & CPU    & Iter                      & CPU       & Iter                      & CPU       & Iter                      & CPU& Iter                      & CPU       & Iter                      & CPU       \\ \hline
			\multicolumn{1}{|c|}{$2^{6} - 1$}   & 3969   &    32                   &  0.19    &   32                       &  0.13     &                    32     &   0.097    &     32                     &   0.10  & 32&0.095 & 32&0.10  \\ \hline
			\multicolumn{1}{|c|}{$2^{8} - 1$}  & 65025  &  32                     & 4.18   &   32                       &  4.15     &  32                        &  4.05      &    32                      &   4.10& 32& 3.85 & 32&3.82   \\ \hline
			\multicolumn{1}{|c|}{$2^{10} - 1$}  & 1046529  &  32                     &  846.45   &  32 &    839.50      &  32                       &  840.90   &    32                      &  841.57 & 32&815.53 & 32&817.22  \\ \hline
			\multicolumn{1}{|c|}{$2^{12} - 1$}  & 16769025 & -                      & $>3600$  & -  & $>3600$           &  -                       & $>3600$   &    -                     &   $>3600$ & -&$>3600$ & -&$>3600$ \\ \hline
			\multicolumn{1}{|c|}{$2^{14} - 1$}  & 268402689 &   -                     & $>3600$  &   -         & $>3600$   &              -            &  $>3600$    &     -                     &  $>3600$ & -&$>3600$ & -&$>3600$  \\ \hline
		\end{tabular}
  }}

\end{subtable}

\begin{subtable}[t]{\textwidth}{
		\caption{GMRES with the PRESB}
\resizebox{\textwidth}{!}{
		\begin{tabular}{|cc|cc|cc|cc|cc|cc|cc|}
			\hline
			\multicolumn{2}{|c|}{$\textrm{tol} = 10^{-8}$} & \multicolumn{2}{c|}{$(\alpha, \beta)=(-600,150)$} & \multicolumn{2}{c|}{$(\alpha, \beta)=(-100,-25)$} & \multicolumn{2}{c|}{$(\alpha, \beta)=(100,-100)$} & \multicolumn{2}{c|}{$(\alpha, \beta)=(-100,100)$}& \multicolumn{2}{c|}{$(\alpha, \beta)=(-100,1)$} & \multicolumn{2}{c|}{$(\alpha, \beta)=(1,-100)$} \\ \hline
			$M$                                 & DoF      & Iter                   & CPU    & Iter                      & CPU       & Iter                      & CPU       & Iter                      & CPU & Iter                      & CPU       & Iter                      & CPU       \\ \hline
			\multicolumn{1}{|c|}{$2^{6} - 1$}   & 7938   &    176                   &  0.50  &    169                      &   0.42    &   175                       &  0.44     &    167                    &    0.40 &  176 & 0.41 & 172 & 0.37     \\ \hline
			\multicolumn{1}{|c|}{$2^{8} - 1$}  & 130050  &    872                   & 39.28   &    840                      &   38.57    &     848                     &   37.66    &     858                     &   39.18 & 1026& 43.66 & 890 &39.03   \\ \hline
			\multicolumn{1}{|c|}{$2^{10} - 1$}  &  2093058  &   -                     & $>3600$  &   -         & $>3600$   &              -            &  $>3600$    &     -                     &  $>3600$ & -&$>3600$ & -&$>3600$  \\ \hline
			\multicolumn{1}{|c|}{$2^{12} - 1$}  & 33538050 &   -                     & $>3600$  &   -         & $>3600$   &              -            &  $>3600$    &     -                     &  $>3600$ & -&$>3600$ & -&$>3600$  \\ \hline
			\multicolumn{1}{|c|}{$2^{14} - 1$}  & 536805378 &   -                     & $>3600$  &   -         & $>3600$   &              -            &  $>3600$    &     -                     &  $>3600$ & -&$>3600$ & -&$>3600$ \\ \hline
		\end{tabular}
  }}
	
		\end{subtable}

\begin{subtable}[t]{\textwidth}{
		\caption{GMRES with the PMHSS}
\resizebox{\textwidth}{!}{
		\begin{tabular}{|cc|cc|cc|cc|cc|cc|cc|}
			\hline
			\multicolumn{2}{|c|}{$\textrm{tol} = 10^{-8}$} & \multicolumn{2}{c|}{$(\alpha, \beta)=(-600,150)$} & \multicolumn{2}{c|}{$(\alpha, \beta)=(-100,-25)$} & \multicolumn{2}{c|}{$(\alpha, \beta)=(100,-100)$} & \multicolumn{2}{c|}{$(\alpha, \beta)=(-100,100)$}& \multicolumn{2}{c|}{$(\alpha, \beta)=(-100,1)$} & \multicolumn{2}{c|}{$(\alpha, \beta)=(1,-100)$} \\ \hline
			$M$                                 & DoF      & Iter                   & CPU    & Iter                      & CPU       & Iter                      & CPU       & Iter                      & CPU& Iter                      & CPU       & Iter                      & CPU       \\ \hline
			\multicolumn{1}{|c|}{$2^{6} - 1$}   & 3969   &    33                   &  0.13   &   32                      &  0.11     &                    33     &   0.13    &     32                     &   0.13  & 32 &0.43 & 32 &0.44  \\ \hline
			\multicolumn{1}{|c|}{$2^{8} - 1$}  & 65025  &  72                     & 13.13  &   71                       &  13.06     &  70                        &  13.11      &    70                     &   13.13 & 71& 13.19 & 71&13.24   \\ \hline
			\multicolumn{1}{|c|}{$2^{10} - 1$}  & 1046529   &-                     & $>3600$  &   -         & $>3600$   &              -            &  $>3600$    &     -                     &  $>3600$ & -&$>3600$ & -&$>3600$  \\ \hline
			\multicolumn{1}{|c|}{$2^{12} - 1$}  & 16769025 & -                     & $>3600$  &   -         & $>3600$   &              -            &  $>3600$    &     -                     &  $>3600$ & -&$>3600$ & -&$>3600$  \\ \hline
			\multicolumn{1}{|c|}{$2^{14} - 1$}  & 268402689 &   -                     & $>3600$  &   -         & $>3600$   &              -            &  $>3600$    &     -                     &  $>3600$ & -&$>3600$ & -&$>3600$  \\ \hline
		\end{tabular}
  }}

\end{subtable}

	\caption{Numbers of iterations and CPU times for Example \ref{ex:var}}
	\label{table_2D_exp1_var_A_iter}
\end{table}


\section{Conclusions}\label{sec:conclusions}
In this study, we have devised optimal preconditioners for efficiently solving complex-shifted Laplacian systems, considering both constant and variable coefficient cases. For the constant coefficient case, our approach ensures optimal convergence when used in conjunction with the MINRES method, while maintaining computational efficiency. For the variable coefficient case, we have demonstrated the convergence of our algorithm under the relevant assumptions. To validate the efficacy of our proposed preconditioners, we conducted numerical experiments. The results demonstrate that our preconditioners perform as expected, further confirming their effectiveness.

{Although our preconditioner relies on the uniform grid discretization only existing in regular domain,  we will consider extending our preconditioning technique to the problem on irregular domain via domain decomposition method in the future.}

\section*{Acknowledgments}
The work of Sean Hon was supported in part by the Hong Kong RGC under grant 22300921, a start-up grant from the Croucher Foundation, and a Tier 2 Start-up Grant from Hong Kong Baptist University. 	The work of Xue-lei  Lin was partially supported by research grants: 2021M702281 from China Postdoctoral Science Foundation, 12301480 from NSFC,  HA45001143 from Harbin Institute of Technology, Shenzhen, HA11409084  from Shenzhen.

\section*{Data Availability Statement}
The data of natural images and codes involved in this paper are available from the corresponding author on reasonable request.

\section*{Declarations}
\begin{description}
	\item[1.]I confirm that I have read, understand and agreed to the submission guidelines, policies and submission 
	declaration of the journal.\\
	\item[2.]I confirm that all authors of the manuscript have no conflict of interests to declare.\\
	\item[3.] I confirm that the manuscript is the authors' original work and the manuscript has not received prior publication 
	and is not under consideration for publication elsewhere.\\
	\item[4.] On behalf of all Co-Authors, I shall bear full responsibility for the submission.\\
	\item[5.] I confirm that all authors listed on the title page have contributed significantly to the work, have read the 
	manuscript, attest to the validity and legitimacy of the data and its interpretation, and agree to its submission.\\
	\item[6.]  I confirm that the paper now submitted is not copied or plagiarized version of some other published work.\\
	\item[7.] I declare that I shall not submit the paper for publication in any other Journal or Magazine till the decision is 
	made by journal editors. \\
	\item[8.] If the paper is finally accepted by the journal for publication, I confirm that I will either publish the paper
	immediately or withdraw it according to withdrawal policies.\\
	\item[9.] I understand that submission of false or incorrect information/undertaking would invite appropriate penal 
	actions as per norms/rules of the journal.
\end{description}
 \nocite{*}
 
\bibliographystyle{plain}      
\bibliography{ref}

\end{document}